\documentclass[12pt,a4paper]{amsart}

\usepackage{amssymb,amsmath,amsthm}
\usepackage{color}
\usepackage{rotating}
\usepackage{fullpage}
\usepackage{comment}
\usepackage{hyperref}
\usepackage{soul}  
\usepackage[all]{xy}
\usepackage{tikz}
\usepackage{tikz-cd}
\DeclareMathOperator{\inte}{int}

\newtheorem{theorem}{Theorem}
\newtheorem{lemma}[theorem]{Lemma}

\theoremstyle{definition}
\newtheorem{definition}[theorem]{Definition}
\newtheorem{remark}[theorem]{Remark}

\title{From the Sharkovskii theorem to periodic orbits for the R\"ossler system}

\author{Anna Gierzkiewicz}
\email{anna.gierzkiewicz@urk.edu.pl}
\address{Department of Applied Mathematics, University of Agriculture in Krak\'ow,
ul. Balicka 253c,
30--198 Krak\'ow, Poland
}

\author{Piotr Zgliczy\'nski}
\email{umzglicz@cyf-kr.edu.pl}
\address{Institute of Computer Science, Jagiellonian University,
ul. \L ojasiewicza 6, 30-348 Krak\'ow, Poland
}

\begin{document}

\begin{abstract}

We extend Sharkovskii's theorem to the cases of $N$-dimensional maps which are close to 1D maps, with an attracting
$n$-periodic orbit. We prove that, with relatively weak topological assumptions, there exist also $m$-periodic orbits for all $m\triangleright n$ in Sharkovskii's order, in the nearby.

We also show, as an example of application, how to obtain such a result for the R\"ossler system with an attracting periodic orbit, for four sets of parameter values. The proofs are computer-assisted.


\end{abstract}


\maketitle

\section{Introduction}

Let us begin by recalling the  Sharkovskii Theorem \cite{ShU,Stefan,Block}:

\begin{theorem}[Sharkovskii]\label{th:shar}
Define an ordering `$\triangleleft$' of natural numbers:
\begin{equation}\label{eq:order}
\begin{array}{r@{\ \triangleleft\ }c@{\ \triangleleft\ }c@{\ \triangleleft\ }c@{\ \triangleleft\ }l}
3\triangleleft 5 \triangleleft 7 \triangleleft 9 \triangleleft\ \dots & 2\cdot 3 & 2 \cdot 5 & 2 \cdot 7 & \dots
\\
 \dots & 2^2\cdot 3 & 2^2 \cdot 5 & 2^2 \cdot 7 & \dots
\\
 \dots & 2^3\cdot 3 & 2^3 \cdot 5 & 2^3 \cdot 7 & \dots
\\
 \multicolumn{5}{c}{\dots \dots \dots \dots \dots \dots \dots \dots \dots \dots }
\\
 \dots & 2^{k+1} & 2^k & 2^{k-1} & \ldots\ \triangleleft 2^2 \triangleleft 2 \triangleleft 1.
\\
\end{array}
\end{equation}
Let $f: I \to \mathbb{R}$ be a continuous map of an interval. If $f$ has an $n$-periodic point and $n\triangleleft m$, then $f$ also has an $m$-periodic point.
\end{theorem}

The paper addresses the following question: Can the Sharkovskii theorem  be carried over to multidimensional dynamics? Without additional assumptions the answer is negative.
However, if the map is in some sense close to 1-dimensional then one can hope for a positive answer.

The typical situation in which such a result can be of interest is the following: consider a Poincar\'e map $F$ for some ODE with the following
properties: there exists an attracting set which is close to being one-dimensional. Then we have an approximate 1D map $f$ which may posses some periodic
orbits to which the Sharkovskii theorem applies, implying the existence  of some periodic points for $f$. We would like to infer that $F$ has
also orbits with the same periods. The result of this type has been given in \cite{PZszarI, PZmulti}, but it can hardly
be directly applied to an explicit example (the R\"ossler system considered in our paper) as it is difficult to obtain any reasonable quantitative condition
on the difference between $F$ and $f$ from the proof in \cite{PZszarI}.

However, in our recent paper \cite{AGPZrossler}, where we consider two sets of parameters in the R\"ossler system with attracting periodic orbits: 3-periodic for the first set and 5-periodic for the other, we proved the existence of the periods implied by the Sharkovskii theorem quite easily.
Just as in \cite{PZszarI}, the proof relies on the construction of a set of covering relations, but due to the fact that orbit we started with was attracting the construction was much simpler in this case.

The goal of this work the generalize this observation, which is the main abstract result in this paper.
\begin{theorem}
\label{thm:sh-manyD}
Consider a continuous map $F:I\times \overline{B}(0,R)\to \inte \left(I\times \overline{B}(0,R)\right)$, where $I\subset \mathbb{R}$ is a closed interval and $\overline{B}(0,R) \subset \mathbb{R}^{n-1}$ a closed ball of radius $R$. Let us denote by $(x,y)$
points in $I\times \overline{B}(0,R)$.

 Suppose that $F$ has an $n$-periodic point $(x_0,y_0) \in \mathbb{R} \times \mathbb{R}^{n-1}$ with least period $n$ and denote its orbit by $\{(x_0,y_0), (x_1,y_1)=F(x_0,y_0), \dots, (x_{n-1},y_{n-1}) = F^{n-1}(x_0,y_0), (x_n,y_n)=(x_0,y_0)\}\subset \inte I\times \overline{B}(0,R)$.

Suppose that there exist $\delta_0$, $\delta_1$, \dots, $\delta_{n-1}>0$ such that
\[
\forall\: i\in\{0,\dots,n-1\} \qquad F\left([x_i\pm\delta_i]\times\overline{B}(0,R)\right) \subset (x_{i+1}\pm\delta_{i+1})\times B(0,R).
\]
Then for every natural number $m$ succeeding $n$ in the Sharkovskii order \eqref{eq:order} $F$ has a point with the least period $m$.
\end{theorem}
The geometric situation in which the above theorem is applicable is as follows:
 consider a map $F: \mathbb{R}^n \to \mathbb{R}^n$, assume that $F$ has an nearly one-dimensional attractor  $\mathcal{A}$, by which we mean that in suitable coordinates  $\mathcal{A}$ is contained in $I \times \overline{B}(0,R)$ for some small $R>0$ and $F\left(I\times \overline{B}(0,R)\right)\subset \inte \left(I\times \overline{B}(0,R)\right)$, and $F$ has an attracting orbit of period $m$ in $\inte \left(I\times \overline{B}(0,R)\right)$, which basin of attraction contains $C_i$ products of interval and a ball $[x_i -\delta_i,x_i+\delta] \times \overline{B}(0,R)$ such that $F(C_i) \subset C_{i+1}$ for $i=0,\dots,n-1$.

The essential difference between the construction used here and in \cite{PZszarI} is the version of the proof of the Sharkovskii theorem
which is used to construct some covering relations. In the present work we use the ideas from the proof of Burns and Hasselblatt \cite{Burns}
while in \cite{PZszarI} is was based the so-called \u{S}tefan cycle \cite{Stefan, Block}.

As an application of Theorem~\ref{thm:sh-manyD} we study  the R\"ossler system \cite{AGPZrossler} for different values of parameters  with attracting periodic orbits of  $3$, $5$ and $6$.  We  verify the assumptions of Theorem~\ref{thm:sh-manyD}, hence we obtain an infinite number of periodic orbits. For precise statements see Sec~\ref{sec:App}). The proofs for the R\"ossler system are computer-assisted, written in C++ with the use of CAPD (Computer-Assisted Proofs in Dynamics) library \cite{capd,capd-article} for interval arithmetic, differentiation and ODE integration.

The content of the paper can be described as follows. In Section \ref{sec:BHproof} we recall some ideas and facts from the proof of the Sharkovskii theorem by   Burns and Hasselblatt \cite{Burns}. In Section~\ref{sec:coveringRn} we discuss the notion of covering relations, which is the main tool used in this work to obtain periodic orbits. In Section~\ref{sec:grids} we define a notion of a \emph{contracting grid} around a periodic orbit.   In Section~\ref{sec:mainThm}  we  prove the main theorem. Finally in Section~\ref{sec:App} we apply our result to the R\"ossler system from \cite{AGPZrossler}.

\subsection*{Notation:}
\begin{itemize}
	\item We use the common notation for the closure, interior, and boundary of a topological set $A\subset\mathbb{R}^k$, which are  $\overline{A}$, $\inte A$, and $\partial A$, respectively.
	\item $\pi_i$ denotes the projection onto the $i$th coordinate in $\mathbb{R}^N$.
	\item By an `$n$-periodic orbit' or `point', we understand an orbit or a point with basic period $n$.
\end{itemize}

\section{Proof of Sharkovskii's theorem by Burns \& Hasselblatt}\label{sec:BHproof}

First we introduce some necessary definitions taken directly from \cite{Burns}.

\subsection{Interval covering relation `$\longrightarrow$'}

Let $\mathcal{O} \subset \mathbb{R}$ be a set of $n$ points.
By an \emph{$\mathcal{O}$-interval} we understand any interval $J\subset [\min\mathcal{O}, \max\mathcal{O}]$ of positive length with the endpoints in $\mathcal{O}$, that is, $\partial J \subset \mathcal{O}$.

Let $f:\mathcal{O} \to \mathcal{O}$ be an $n$-periodic permutation with every $x\in \mathcal{O}$ being an $n$-periodic point for $f$.

In the paper we will use a stronger notion of a one-dimensional covering relation between intervals than the classical one used by Block \textit{et al.} \cite{Block}. Our definition appears in the paper of Burns and Hasselblatt \cite{Burns} as the \emph{$\mathcal{O}$-forced covering} between $\mathcal{O}$-intervals, which we simply call the covering between $\mathcal{O}$-intervals:

\begin{definition}[($\mathcal{O}$-forced) covering relation]\label{def:Oforced_covering}
An $\mathcal{O}$-interval $I$ \emph{covers} an $\mathcal{O}$-interval $J$ (denoted by $I \overset{f}{\rightarrow} J $ or simply  $I \rightarrow J $), if
there exists a $\mathcal{O}$-subinterval $K\subset I$ such that
\begin{equation}
\min f(\partial K) \leq \min J \text{ \qquad and \qquad } \max J \leq \max f(\partial K).
\end{equation}
\end{definition}

The above relation fulfills the Itinerary Lemma \cite{Block,Burns}, which is crucial in proving the existence of periodic orbits in the proofs of Sharkovskii's theorem:
\begin{theorem}[Itinerary Lemma]
\label{th:1d-covering}
Let $f: I \to \mathbb{R}$ be a continuous map on an interval $I \subset\mathbb{R}$ and $\mathcal{O}\subset I$ be an $n$-periodic point for $f$. Assume that we have a sequence of $\mathcal{O}$-intervals $J_j \subset I$ for $j=0,\dots,m-1$ such that
\begin{equation}\label{eq:1d_loop}
  J_0  \overset{f}{\longrightarrow} J_1 \overset{f}{\longrightarrow} J_2 \overset{f}{\longrightarrow} \dots \overset{f}{\longrightarrow} J_{m-1} \overset{f}{\longrightarrow} J_0.
\end{equation}
Then there exists a point $x \in J_0$, such that $f^j(x) \in J_j$ for $j=1,\dots,m-1$ and $f^m(x)=x$.
\end{theorem}
A point which fulfills the thesis of Theorem \ref{th:1d-covering} is said to \emph{follow the loop} \eqref{eq:1d_loop}. A loop of length $m$ (such as \eqref{eq:1d_loop}), we will call shortly an $m$-loop.

In general, from Theorem \ref{th:1d-covering} we do not know whether the point's period is fundamental. To make sure that the point $x\in J_0$ following the loop \eqref{eq:1d_loop} is indeed $m$-periodic, we must add some assumptions on the $\mathcal{O}$-intervals forming the loop of covering relations. Particular criteria may put additional assumptions on the intervals $J_i$ or some restrictions on their order. We use the following criterion, which includes all the cases studied in \cite{Burns}:

\begin{definition}[\cite{Burns}, Lemma 2.6\footnote{It is a more specific definition than the one of an \emph{elementary loop} from \cite[Def. 2.4]{Burns}.}]\label{def:basic_loop}
We call a loop \eqref{eq:1d_loop} a \emph{non-repeating loop}, if it fulfills the following conditions:
\begin{enumerate}
	\item it is not followed by any endpoint $x\in\bigcup_{i=0}^{m-1} \partial J_i$;
	\item $\displaystyle \inte J_0 \cap \bigcup_{i=1}^{m-1} J_{i} = \varnothing$.
\end{enumerate}
\end{definition}

\begin{lemma}[\cite{Burns}, Lemma 2.6]\label{lem:non-repeating}
If a loop \eqref{eq:1d_loop} is non-repeating, then every point following it has period $m$.\qed
\end{lemma}

\subsection{Proposition 6.1 of \cite{Burns}}

Let now $f: I \to \mathbb{R}$ be a continuous map on an interval $I \subset\mathbb{R}$ and $\O\subset I$ be an $n$-periodic orbit for $f$.

The following theorem is Proposition 6.1 of \cite{Burns} for the map $f$.
\begin{theorem}[\cite{Burns}]\label{th:Burns}
	For every number $m$ succeeding $n$ in Sharkovskii's order \eqref{eq:order} (\textit{i.e.} $n\triangleleft m$) there exists a non-repeating $\mathcal{O}$-forced $m$-loop of $\mathcal{O}$-intervals $J_i$, $i=0,\dots, m-1$:
	\begin{equation}
	  J_0  \overset{f}{\longrightarrow} J_1 \overset{f}{\longrightarrow} J_2 \overset{f}{\longrightarrow} \dots \overset{f}{\longrightarrow} J_{m-1} \overset{f}{\longrightarrow} J_0\text{,}
	\end{equation}
	which proves the existence of an $m$-periodic point for $f$ in the interval $ [\min \mathcal{O},\max\mathcal{O}]$.
\end{theorem}

Since the above theorem is crucial for our construction, we present here the outline of the proof.

\begin{remark}[Sketch of the proof]
The proof of above theorem, which is the main result of \cite{Burns}, relies by induction on either:
\begin{itemize}
	\item constructing the so-called \u{S}tefan sequence of some length $2\leq l\leq n$ of $\mathcal{O}$-intervals $J_i$, $i=0,\dots,l-1$, which form the diagram of $\mathcal{O}$-forced covering relations ($\inte J_0$ is disjoint with all other $J_i$'s):
	\vspace{3mm}
	\begin{equation}\label{eq:diagOdd}
\begin{tikzcd}
	&
	J_1
		\arrow[r, dashed]
		\arrow[loop, distance=2.5em, in=100, out=170]
	&
	\text{\rotatebox[origin=c]{-20}{$\ldots$}}
		\arrow[rd, dashed]
	&
\\
	J_0
		\arrow[d, shift left]
		\arrow[ru]
		\arrow[rru, "\text{\tiny to $J_{l-2k-1}$}"', dashed]
		\arrow[rrr] \arrow[rrdd]
	& & 	&
	J_{l-5}
		\arrow[d]
\\
	J_{l-1}
		\arrow[u, shift left]
	& &	&
	J_{l-4} 	
		\arrow[ld]
\\
    &
    J_{l-2}
    		\arrow[lu]
    	&
    	J_{l-3}
    		\arrow[l]
    	&
\end{tikzcd}
\end{equation}	
From the above diagram \eqref{eq:diagOdd} one deduces the existence of $m$-periodic points, from non-repeating loops:
\begin{itemize}
	\item $m=1$ from \qquad $
	  J_1 \longrightarrow J_1$,
	\item $m \geq l$ from \qquad $
	  J_0 \rightarrow J_1 \underbrace{\dots \rightarrow J_1}_{\text{$(m-l) \times$ `$\rightarrow J_1$'}}  \rightarrow J_2 \rightarrow \dots \rightarrow J_{l-1} \rightarrow J_0$,
	\item even $m<l$ from \qquad $
	  J_0 \rightarrow  J_{l-m+1} \rightarrow  J_{l-m+2} \rightarrow \dots \rightarrow J_{l-1} \rightarrow J_0$.
\end{itemize}
	\item or reducing to the case of the orbit of length $n/2$ for the map $f^2$, if the construction of a \u{S}tefan sequence is impossible.
\end{itemize}
\end{remark}

Although the authors formulate their Proposition 6.1 in a more general way, the loops which appear in the proof of Theorem \ref{th:Burns} are in fact non-repeating and the proof relies on that fact. Let us review shortly the induction on the length of $\mathcal{O}$, $n\in\mathbb{N}$ presented in the proof of \cite[Proposition 6.1]{Burns}:

\begin{enumerate}
	\item For $n=1$ there are no loops of length $m\triangleright n$.
	\item For $n=2$ the constructed loop is of the form
	\[
	J_1 \rightarrow J_1\text{,}
	\]
	which is non-repeating.
	\item Suppose now that for all lengths of $\mathcal{O}$ smaller that $n$ the loops constructed in the proof are non-repeating. For the length $n>2$ and a number $m\triangleright n$ there are two possible cases.
	\begin{enumerate}
		\item We are able to construct a \u{S}tefan sequence and the non-repeating $m$-loop is of the form:
		\[
	  J_0 \rightarrow \underbrace{\dots \rightarrow J_1}_{\geq 0 \text{ times}} \rightarrow J_1 \rightarrow J_2 \rightarrow \dots \rightarrow J_{k} \rightarrow J_0\text{,}
	  \]
	  with $\inte J_0 \cap \bigcup_{i=1}^{k} J_i = \varnothing$ and $k\geq 1$.
	  \item We cannot construct a \u{S}tefan sequence, but then either $m=1$ and the case is trivial, or $n$, $m$ are even. Hence, by induction, we have a non-repeating loop of length $\frac{m}{2}$ for the map $f^2$:
	  \begin{equation}\label{eq:even_loop}
	  J_0 \overset{f^2}\longrightarrow  J_{1} \overset{f^2}\longrightarrow  J_{2} \overset{f^2}\longrightarrow \dots \overset{f^2}\longrightarrow J_{\frac{m}2-1} \overset{f^2}\longrightarrow J_0.
	  \end{equation}
	  Next, the construction from the proof of \cite[Proposition 6.1]{Burns} extends the above $\frac{m}2$-loop \eqref{eq:even_loop} to the following $m$-loop for $f$:
	  \begin{equation*}
	  J_0
	  \overset{f}\longrightarrow J'_{0}
	  \overset{f}\longrightarrow J_{1}
	  \overset{f}\longrightarrow J'_{1}
	  \overset{f}\longrightarrow
	  \dots
	  \overset{f}\longrightarrow J_{\frac{m}2-1}
	  \overset{f}\longrightarrow J'_{\frac{m}2-1}
	  \overset{f}\longrightarrow J_0\text{,}	
	  \end{equation*}
	  where $\bigcup_{i=0}^{m/2} J_i \cap \bigcup_{i=0}^{m/2} J'_i = \varnothing$. In particular, $J_0 \cap J'_{i} = \varnothing$ for $i=0,\dots,\frac{m}2-1$.
	\end{enumerate}
\end{enumerate}

\subsection{Proper covering}

Let us introduce now a more specific notion of covering between $\mathcal{O}$-intervals, which we can easily compare later to the notion of horizontal covering between segments (Section \ref{sec:coveringRn}).

\begin{definition}[Proper covering of intervals]
Assume that $I, J \subset \mathbb{R}$ are $\mathcal{O}$-intervals. We say that
$I$ \emph{$f$-covers $J$ properly} (denoted\footnote{The same symbol is used in \cite{Burns}, but for a different notion.} by $I \overset{f}{\rightarrowtail} J$ or simply  $I \rightarrowtail J $), if 
\begin{equation}\label{eq:proper_cover}
\min f(\partial I) \leq \min J \text{ \qquad and \qquad } \max J \leq \max f(\partial I).
\end{equation}
\end{definition}

It is easy to see from the definition of $\mathcal{O}$-forced covering (Definition \ref{def:Oforced_covering}) that the following statement is true.

\begin{lemma}\label{lem:proper_covering}
If $I\rightarrow J$ is an $\mathcal{O}$-forced covering relation, then there exists an $\mathcal{O}$-subinterval $K\subset I$ such that $K\rightarrowtail J$.\qed
\end{lemma}

We can finally replace $\mathcal{O}$-forced loops by loops of proper covering relations.

\begin{lemma}\label{lem:proper_loop}
	For every number $m$ succeeding $n$ in Sharkovskii's order \eqref{eq:order} (\textit{i.e.} $n\triangleleft m$) there exists a non-repeating $m$-loop of proper covering relations between $\mathcal{O}$-intervals $K_i$, $i=0,\dots, m-1$:
	\begin{equation}\label{eq:proper_loop}
	  K_0  \overset{f}{\rightarrowtail} K_1 \overset{f}{\rightarrowtail} K_2 \overset{f}{\rightarrowtail} \dots \overset{f}{\rightarrowtail} K_{m-1} \overset{f}{\rightarrowtail} K_0\text{.}
	\end{equation}
\end{lemma}

\begin{proof}
Fix $m\triangleright n$ and consider the non-repeating $\mathcal{O}$-forced loop obtained from Theorem \ref{th:Burns}:
	\begin{equation*}
	  J_0  \overset{f}{\longrightarrow} J_1 \overset{f}{\longrightarrow} J_2 \overset{f}{\longrightarrow} \dots \overset{f}{\longrightarrow} J_{m-1} \overset{f}{\longrightarrow} J_m=J_0\text{.}
	\end{equation*}
	Now apply Lemma \ref{lem:proper_covering} to every relation $J_i \rightarrow J_{i+1}$, $i=0,\dots, m-1$. We obtain $m$ $\mathcal{O}$-intervals $K_i\subset J_i$, for $i=0,\dots, m-1$ which fulfill the proper covering relations:
	\[
	K_i \rightarrowtail J_{i+1}\text{, } \quad i=0,\dots, m-1.
	\]
	Observe now that also each of the proper covering relations  $K_i \rightarrowtail K_{i+1}$ are true, because each $K_{i+1}$ is a subinterval of $J_{i+1}$.
	Finally, note that the loop 			
	\begin{equation*}
	  K_0  \overset{f}{\rightarrowtail} K_1 \overset{f}{\rightarrowtail} K_2 \overset{f}{\rightarrowtail} \dots \overset{f}{\rightarrowtail} K_{m-1} \overset{f}{\rightarrowtail} K_0\text{}
	\end{equation*}
	is non-repeating, because both conditions (1) and (2) from Definition \ref{def:basic_loop} are fulfilled if we replace the $\mathcal{O}$-intervals by their $\mathcal{O}$-subintervals.
\end{proof}

\section{Horizontal covering relation `$\Longrightarrow$'}\label{sec:coveringRn}

 In order to have the Itinerary Lemma in many dimensions we need a good notion of covering where we have a direction of possible expansion
and an apparent contraction in other directions.

For this end  we recall  the notion of covering for h-sets in $\mathbb{R}^N$ from \cite{PZszarI,PZmulti,AGPZrossler}.
It is, in fact, a particular case of the similar notion from \cite{GZ}, but with exactly one exit (or 'unstable')  direction. For such a covering relation a special version of Itinerary Lemma is true (Th.\ \ref{th:periodic} below) and we use it to prove the existence of periodic points for multidimensional maps.

\begin{definition}
An \emph{h-set} is a hyper-cuboid $S=[a_1,b_1]\times \dots \times [a_N,b_N]\subset \mathbb{R}^N$ for some $a_i<b_i$, $i=1,\dots,N$, with the following elements distinguished:
\begin{itemize}
	\item its left face \quad $L(S) = \{x\in\partial S \,:\,x_1=a_1\}$,
	\item its right face \quad $R(S) = \{x\in\partial S \,:\,x_1=b_1\}$,
	\item its horizontal boundary \quad $H(S) = \overline{\partial S \setminus (L(S)\cup R(S))}$,
	\item its left side \quad $\mathcal{L}(S) = \{x\in\mathbb{R}^N \,:\,x_1<a_1\}$,
	\item its right side \quad $\mathcal{R}(S) = \{x\in\mathbb{R}^N \,:\,x_1>b_1\}$.
\end{itemize}
\end{definition}

Now we define a horizontal covering relation between h-sets.

\begin{figure}[h]
	\includegraphics[height=7cm]{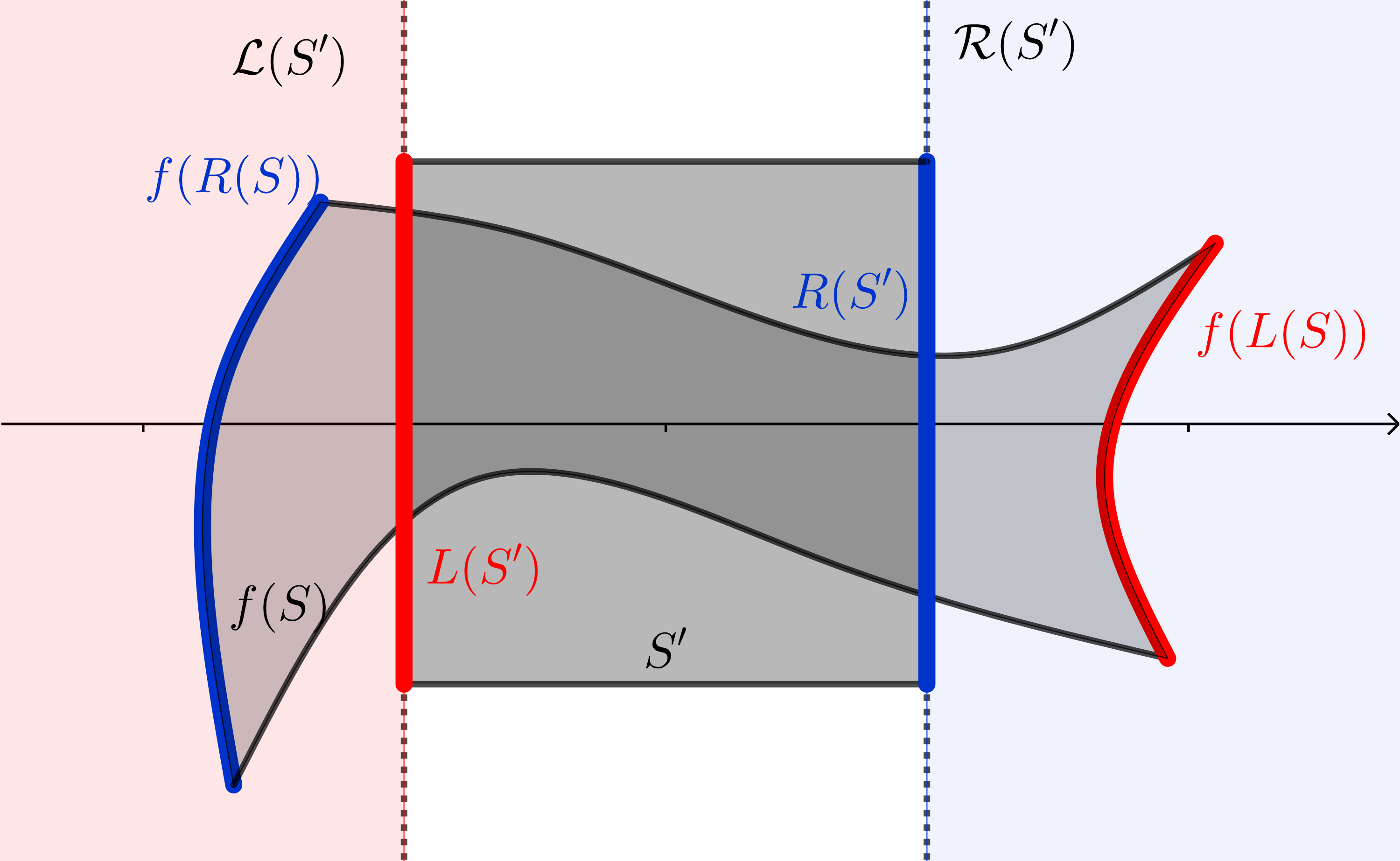}
	\caption{\label{fig:cover}Horizontal covering between 2D h-sets $S \overset{f}{\Longrightarrow} S'$}
\end{figure}

\begin{figure}[h]
	\includegraphics[width=0.7\textwidth]{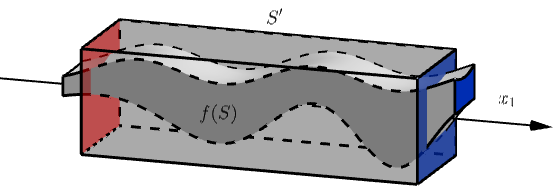}
	\caption{\label{fig:cover3d}Horizontal covering between 3D h-sets $S \overset{f}{\Longrightarrow} S'$}
\end{figure}

\begin{definition}
\label{def:cov}
Let $S$, $S'$ be two h-sets and $f: V \to \mathbb{R}^N$ be a continuous map on an open neighborhood of $S\subset V\subset\mathbb{R}^N$.
We say that $S$ \emph{$f$-covers $S'$ horizontally} and denote by $S \overset{f}{\Longrightarrow} S'$ if
\begin{equation}\label{eq:cond9}
f(S) \subset \left(\mathcal{L}(S')\cup S'\cup \mathcal{R}(S')\right) \setminus H(S')\text{,}
\end{equation}
and one of the two conditions hold:
\begin{equation}\label{eq:cond10}
\begin{aligned}
&\text{either } & f(L(S))\subset \mathcal{L}(S') &\text{\quad and \quad} f(R(S))\subset \mathcal{R}(S')\text{,}
\\
&\text{or } & f(L(S))\subset \mathcal{R}(S') &\text{\quad and \quad} f(R(S))\subset \mathcal{L}(S').
\end{aligned}
\end{equation}
\end{definition}
See Fig.\ \ref{fig:cover} for the illustration of horizontal covering in $\mathbb{R}^2$ and Fig.\ \ref{fig:cover3d} for covering in $\mathbb{R}^3$.

Let us emphasize that the above conditions can be easily checked with the use of computer via interval arithmetic and `$<$', `$>$' relations.

The following theorem might be understood as a version of Itinerary Lemma (Theorem~\ref{th:1d-covering}).

\begin{theorem}[\cite{PZszarI}]\label{th:periodic}
Suppose that there occurs a loop of $m$ horizontal $f$-coverings between h-sets $S_i\subset\mathbb{R}^N$, $i=0,\dots,m-1$:
\[
S_0 \overset{f}{\Longrightarrow} S_1\overset{f}{\Longrightarrow} \dots \overset{f}{\Longrightarrow} S_{m-1} \overset{f}{\Longrightarrow} S_m = S_0\text{,}
\]
then there exists $x\in \inte S_0$ such that $f^m(x)=x$ and
\[
\text{for } i=1,\dots ,m-1: \qquad f^i(x)\in \inte S_i.
\]
\end{theorem}

\section{Grids}
\label{sec:grids}

The goal of this section is to introduce the notion of \emph{contracting grid}. In the context of Theorem~\ref{thm:sh-manyD} a contracting grid is a set of cubes of the form $C_i=[x_i-\delta_i,x_i+\delta_i] \times \overline{B}(0,R)$ for $i=0,\dots,n-1$ and the h-sets $S_{ij}=[x_i+\delta_i,x_{j}-\delta_j] \times \overline{B}(0,R)$ for $i,j$, such that $x_i < x_j$, lying between $C_i$'s. Due to the fact that for each $i$ $F(C_i) \subset \inte C_{i+1}$, we obtain, using the ideas from the proof of the Sharkovskii theorem recalled in Section~\ref{sec:BHproof}, a rich set
of horizontal covering relations, which gives us all periods.

\subsection{Model grid}

Let now $\mathcal{O} = \{1, 5,\dots, 4n-3\}$ be the set of $n$ points on the real line $Ox_1$ which, together with $Ox_k$ lines, $k=2,\dots, N$, span the Euclidean space $\mathbb{R}^N = \{x=(x_1, \dots, x_N) \;:\; x_i \in \mathbb{R},\; i=1,\dots, N\}$.
\begin{definition}
By a \emph{model grid} $\mathcal{G}$ enclosing $\mathcal{O}$ we understand an $N$-dimensional full closed hyper-cuboid 
restricted by the following hyper-planes in $\mathbb{R}^N$:
\begin{itemize}
	\item $\{x_1 = 0\}$ and $\{x_1 = 4n-2\}$;
	\item $\{x_k = 1\}$ and $\{x_k = -1\}$, for $k= 2,\dots,N$,
\end{itemize}
with the following elements distinguished:
\begin{itemize}
	\item \emph{inner cubes}: $C_k$, $k\in\mathcal{O}$: \hfil $C_k = \{x\in\mathcal{G} : k-1\leq x_1\leq k+1\}$.
	\\Note that each inner cube $C_k$ contains a single point from $\mathcal{O}$, that is, the point $k$.
	\item \emph{outer segments}: $S_{ij}$,  $i$, $j\in\mathcal{O}$, $i< j$: \hfil $S_{ij} = \{x\in\mathcal{G} : i+1\leq x_1\leq j-1\}$.
	Each outer segment has on its boundary two vertical faces common with two inner cubes $C_i$ and  $C_j$. We will also say that the segment $S_{ij}$ \emph{lies between the cubes} $C_i$ and  $C_j$ or \emph{lies between points} $i$, $j\in\mathcal{O}$.
	\\
	Note also that there is a natural one-to-one correspondence between outer segments and $\mathcal{O}$-intervals on the real line $Ox_1$. The segment $S_{ij}$ corresponding to an $\mathcal{O}$-interval $I_{ij}=[i,j]$ we will denote by $S(I_{ij}) = S_{ij}$.
\end{itemize}
\end{definition}

\begin{figure}[h]
	\includegraphics[height=5cm]{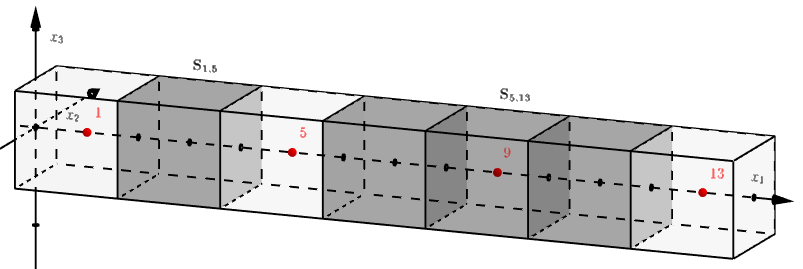}
	\caption{\label{fig:model_grid}A model grid enclosing four points (red) in $\mathbb{R}^3$.
	Two segments are marked in gray: $S_{1,5}$ lying between $1$ and $5$, and $S_{5,13}$ between $5$, $13$. White cubes are the sets $C_1$, $C_5$ and $C_{13}$.}
\end{figure}

For the illustration of the notion of the model grid see Fig. \ref{fig:model_grid}.

\begin{remark}
Note that we can naturally provide a structure of an h-set for a model segment $S=S_{kl}$ lying between the points $k < l \in \{1,5,\dots,4n-3\}$:
\begin{itemize}
	\item $L(S) = \{x\in\partial S \,:\,x_1=k+1\}$,\quad $R(S) = \{x\in\partial S \,:\,x_1=l-1\}$,
	\item $H(S) = \{x\in\partial S : \exists_{i=2,\dots N} |x_i|=1\}$,
	\item $\mathcal{L}(S) = \{x\in\mathbb{R}^N \,:\,x_1<k+1\}$, \quad $\mathcal{R}(S) = \{x\in\mathbb{R}^N \,:\,x_1>l-1\}$.
\end{itemize}
\end{remark}

\subsection{Contracting grid}
Let $U$ be an open subset of $\mathbb{R}^N$.
Consider a set of $n$ points $P = \{p_1, p_2, \dots, p_n\} \subset U$.
\begin{definition}
	A set $G\subset U$ containing $P$ is called a \emph{grid enclosing $P$} if there exists a homeomorphism $h : V\to U$, where $V$ is some open neighborhood of the model grid $\mathcal{G}$ in $\mathbb{R}^N$, such that $P=h(\mathcal{O})$ and $G =h(\mathcal{G})$.
	
	We also define all elements mapped from the model grid: 
	inner cubes and outer segments, as the images of the model elements through $h$.
\end{definition}

Let now $F: U \to \mathbb{R}^N$ be continuous map and $P$ be an $n$-periodic orbit for $F$:
\[
P = \{p_1, p_2, \dots ,p_n, p_{n+1}=p_1\}
\text{, \quad where \quad}
p_1\mapsto p_2\mapsto \dots \mapsto p_n\mapsto p_1.
\]
Suppose that there exists a grid $G = h(\mathcal{G})\subset U$ enclosing the periodic orbit $P$ and denote $F_h = h^{-1}\circ F \circ h : V \to V$.

\begin{definition} We call a grid $G$ a \emph{contracting grid} if
	\begin{equation}\label{eq:contracting}
		F(G) \subset \inte G
		\quad \text{ and } \quad
		\mathop{\forall}_{i=1,\dots,n} F(h(C_i)) \subset \inte h(C_{F_h(i)}) .
	\end{equation}
\end{definition}

\begin{theorem}\label{th:1to2}
	Assume that an $n$-periodic orbit $P$ for $F$ is enclosed by a contracting grid $G= h(\mathcal{G})$ and $\mathcal{O} = h^{-1}(P)$.
	Define $f:\mathcal{O}\to\mathcal{O}$ by $f=F_h|_{\mathcal{O}}$.
	
	If two $\mathcal{O}$-intervals $I$, $I'$ fulfill the proper covering relation $I \overset{f}{\rightarrowtail} I'$, then also the corresponding segments $S(I)$, $S(I')$ fulfill the horizontal covering relation
	\begin{equation}
	S(I) \overset{F_h}{\Longrightarrow} S(I').
	\end{equation}
\end{theorem}

\begin{proof}
Denote by $i<j\in \mathcal{O}$ the endpoints of $I$ and by  $k<l\in \mathcal{O}$ the endpoints of $I'$, so that $I=I_{ij}$, $I'=I_{kl}$, $S(I) = S_{ij}$ and $S(I') = S_{kl}$. From the proper covering $I_{ij} \overset{f}{\rightarrowtail} I_{kl}$ we know that either $f(i)\leq k$ and $f(j)\geq l$,  or $f(i)\geq l$ and $f(j)\leq k$.
\begin{itemize}
	\item \underline{Condition \eqref{eq:cond9}:} The grid is contracting, so $F_h(\mathcal{G})\subset \inte \mathcal{G}$. In particular, $F_h(S(I)) \subset \inte \mathcal{G}\subset\left(\mathcal{L}(S(I'))\cup S(I')\cup \mathcal{R}(S(I'))\right) \setminus H(S(I'))$.
	\item \underline{Condition \eqref{eq:cond10}:}
Consider just the case $f(i)\leq k$ and $f(j)\geq l$, for the other one is analogous.

From the condition \eqref{eq:contracting}, $F_h(C_i)\subset \inte C_{f(i)}$ and  $F_h(C_j)\subset \inte C_{f(j)}$. In particular, $L(S_{ij})\subset C_i$ and $R(S_{ij})\subset C_j$, so $\max \pi_1(F_h(L(S_{ij})))<  f(i) +1 \leq k+1$ and $\min \pi_1(F_h(R(S_{ij})))>  f(j) -1 \geq l-1$, which means that $F_H(L(S_{ij})) \subset \mathcal{L}(S_{kl})$ and  $F_H(R(S_{ij})) \subset \mathcal{R}(S_{kl})$.
\end{itemize}
\end{proof}

\section{The main theorem}
\label{sec:mainThm}

We have the following result.
\begin{theorem}\label{th:Final}
Let $F:U\to \mathbb{R}^N$ be a continuous map on an open set $U\subset\mathbb{R}^N$ and $p_0\in U$ be an $n$-periodic point of $F$ for $n>1$. Suppose that the orbit $O = \{p_0, p_1=F(p_0),\dots p_{n-1} = F^{n-1}(p_0)\}$ is enclosed in a contracting grid $G= h(\mathcal{G})$.

Then for every $m\in\mathbb{N}$ such that $n\triangleleft m$ in Sharkovskii's order \eqref{eq:order}, $F$ has also an $m$-periodic point.
\end{theorem}

\begin{proof}~
Define $f:\mathcal{O}\to\mathcal{O}$ by $f=F_h|_{\mathcal{O}}$, where $\mathcal{O} = h^{-1}(O)$.
From Lemma \ref{lem:proper_loop} we deduce the existence of a non-repeating $m$-loop of proper covering relations between $\mathcal{O}$-intervals $K_i$:	
	\begin{equation}\label{eq:final_loop}
	  K_0  \overset{f}{\rightarrowtail} K_1 \overset{f}{\rightarrowtail} K_2 \overset{f}{\rightarrowtail} \dots \overset{f}{\rightarrowtail} K_{m-1} \overset{f}{\rightarrowtail} K_0\text{.}
	\end{equation}
	
	Next, from Theorem \ref{th:1to2} we deduce that the following loop of horizontal covering relations occurs:
	\begin{equation}\label{eq:Final_Loop}
	  S(K_0)
	  \overset{F_h}{\Longrightarrow}
	  S(K_1)
	  \overset{F_h}{\Longrightarrow}
	  S(K_2)
	  \overset{F_h}{\Longrightarrow}
	  \dots
	  \overset{F_h}{\Longrightarrow}
	  S(K_{m-1})
	  \overset{F_h}{\Longrightarrow}
	  S(K_0)\text{,}
	\end{equation}
	from which follows the existence of a point $x\in S(K_0)$, such that $F_h^m(x)=x$.
	
	The loop \eqref{eq:final_loop} is non-repeating, so one easily observes that the corresponding loop \eqref{eq:Final_Loop} fulfills the condition
\begin{equation}
	S(K_0) \cap \bigcup_{i=1}^{m-1} S(K_i) = \varnothing\text{,}
\end{equation}
	which implies that $m$ is the least period for $x$.
Therefore, the point $h(x)$ is $m$-periodic for $F$.
\end{proof}

\begin{proof}[Proof of Theorem~\ref{thm:sh-manyD}] Observe that under assumptions we have an $n$-periodic orbits which is enclosed
by a contracting grid. The assertion then follows from Theorem~\ref{th:Final}.
\end{proof}


\section{Examples of application: the R\"ossler system}\label{sec:App}

In \cite{AGPZrossler} we proved the existence of $m$-periodic orbits of all $m$ succeeding $n$ in Sharkovskii order \eqref{eq:order} in the R\"ossler system \cite{Rossler76} with an attracting $n$-periodic orbit
\begin{equation}\label{eq:rossler}
\begin{cases}
x'=-y-z,
\\
y'=0.2 y+x,
\\
z'=z (x-a)+0.2
\end{cases}\text{,}
\end{equation}
for two sets of parameters:
\begin{itemize}
	\item $a=5.25$, for which $n=3$ and
	\item $a=4.7$, for which $n=5$.
\end{itemize}
The proof requires finding some explicit special sets (using some ad-hoc trial and error approach) and proving some covering relations between them by computer-assisted methods. However, with Theorem \ref{th:Final}, a part of the proof is easier: the only thing to find and prove its existence is a contracting grid.

\begin{figure}[h]
	\includegraphics[width=0.9\textwidth]{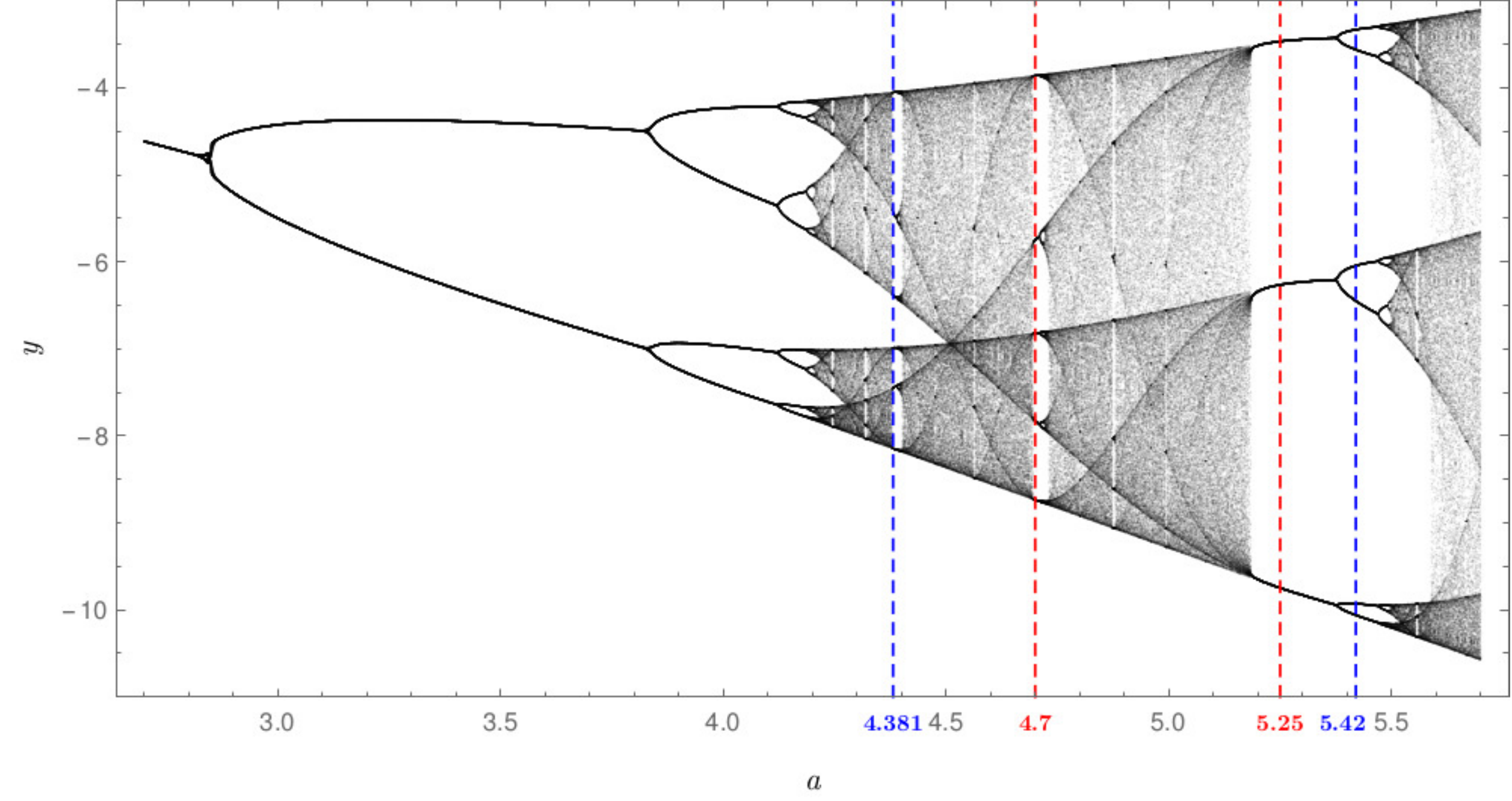}
	\caption{\label{fig:bif}The bifurcation diagram for the Poincar\'e map of the system \eqref{eq:rossler} with $b=0.2$.
	The cases of attracting periodic orbits for $a=4.7$ and $a=5.25$ (red) are treated in \cite{AGPZrossler} and in Subsections \ref{ss:rossler3}, \ref{ss:rossler5}.
	The two other cases of $6$-periodic orbits (blue) are studied in Subsections \ref{ss:rossler6}, \ref{ss:rossler6_2}. }
\end{figure}

We are convinced that it should be possible to find a contracting grid for a wide variety of the parameter $a$'s values for which the system \eqref{eq:rossler} has an attracting periodic orbit.
In this section, we study four values of the parameter $a$: apart from the two cases treated in \cite{AGPZrossler}, we also compare the $6$-periodic orbits appearing in cases $a=4.381$ and $a=5.42$ (see the bifurcation diagram on Fig.\ \ref{fig:bif}). What is interesting, in the latter case one can prove the existence of even more periodic orbits than follow from Theorem \ref{th:Final}. 
This is a reflection of the following fact about forcing relation between periods of interval maps. The set of forced periods depends on the pattern (the permutation
induced on the orbit) and the periodic orbits of the same period may force different sets of periods.  The Sharkovskii Theorem gives only a lower bound
on the set of forced periods. 

In all the four cases, we denote by $\Pi$ the half-plane $\{x=0, y<0\}$ with induced coordinates $(y,z)$, and $P$ is a Poincar\'e map of the system \eqref{eq:rossler} on section $\Pi$, that is the map
\[
P(y,z) = \pi_{(y,z)}\left( \Phi_{T(y,z)}\left(x=0,y,z\right)\right)\text{,}
\]
where $\pi_{(y,z)}$ is the projection on the $(y,z)$ plane, $\Phi_t$ is the dynamical system induced by considered system and $T=T(y,z)$ is a return time, if well-defined.

\subsection{Case $a=5.25$}\label{ss:rossler3}
\begin{figure}[h]
	\includegraphics[width=0.8\textwidth]{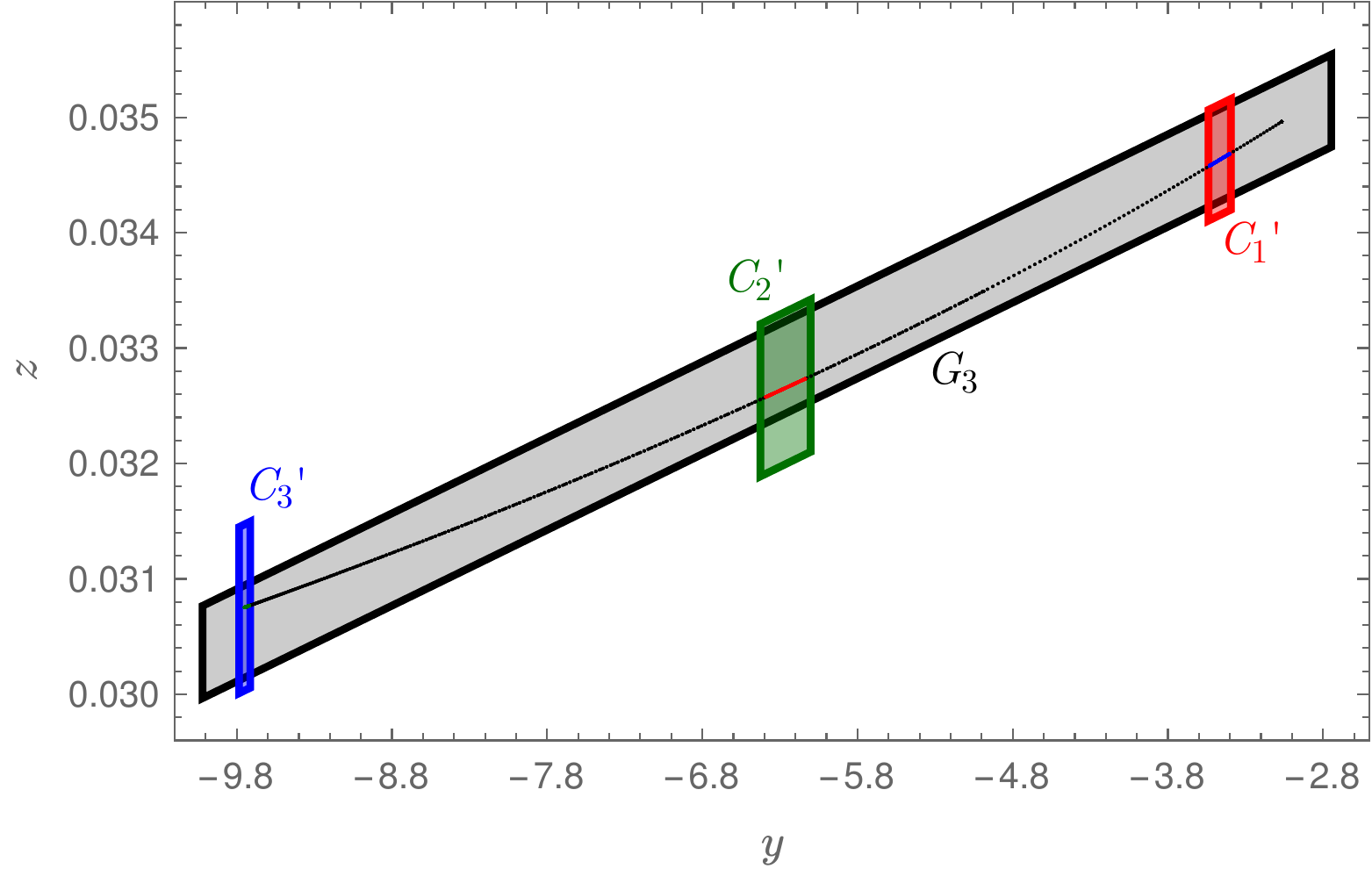}
	\caption{\label{fig:grid3}The contracting grid $G_3$ from Lemma \ref{lem:Roessler525} and its image through $P$.
	\newline The supersets of $C_i$, $i=1,2,3$ and their images are marked in red, green and blue.}
\end{figure}

Let $a=5.25$. Then the system \eqref{eq:rossler} has an attracting $3$-periodic point for $P$ \cite[Lemma 4]{AGPZrossler}.
\begin{lemma}\label{lem:Roessler525}
Let $M=\left[\smallmatrix -1.& 0.000656767 \\ -0.000656767 & -1. \endsmallmatrix\right]$. The parallelogram $G_3$ in the $(y, z)$ coordinates on the section $\Pi$ (see
Fig. \ref{fig:grid3}):
\[
G_3= \bmatrix -6.38401 \\ 0.0327544 \endbmatrix
+ M \cdot
\bmatrix \pm 3.63687 \\ \pm 0.0004 \endbmatrix
\]
is a contracting grid for the attracting $3$-periodic orbit from \cite[Lemma 4]{AGPZrossler}, with the inner cubes $C_i$, $i=1,2,3$, defined by $C_i = C'_i \cap G_3$, where
\begin{align*}
	&C'_1= \bmatrix -3.46642 \\ 0.0346316 \endbmatrix
+ M \cdot
\bmatrix \pm 0.072 \\ \pm 0.00048 \endbmatrix
	\text{,}
	\\
	&C'_2= \bmatrix -6.26401 \\ 0.0326544 \endbmatrix
+ M \cdot
\bmatrix \pm 0.162 \\ \pm 0.00066 \endbmatrix
	\text{,}
	\\
	&C'_3= \bmatrix -9.74889\\ 0.0307529 \endbmatrix
+ M \cdot
\bmatrix \pm 0.036 \\ \pm 0.00072 \endbmatrix.
\end{align*}
\end{lemma}

\begin{proof}
Computer-assisted, \cite{proof}.
\end{proof}

From the above Lemma and Theorem~\ref{th:Final} we obtain
\begin{theorem}(Compare to \cite[Theorem 5]{AGPZrossler})\label{th:r3}
The R\"ossler system \eqref{eq:rossler} with $a=5.25$ has $m$-periodic orbits for any $m\in \mathbb{N}$ in the set $G_3$ defined in Lemma \ref{lem:Roessler525}.
\end{theorem}

\subsection{Case $a=4.7$}\label{ss:rossler5}
\begin{figure}[h]
	\includegraphics[width=0.8\textwidth]{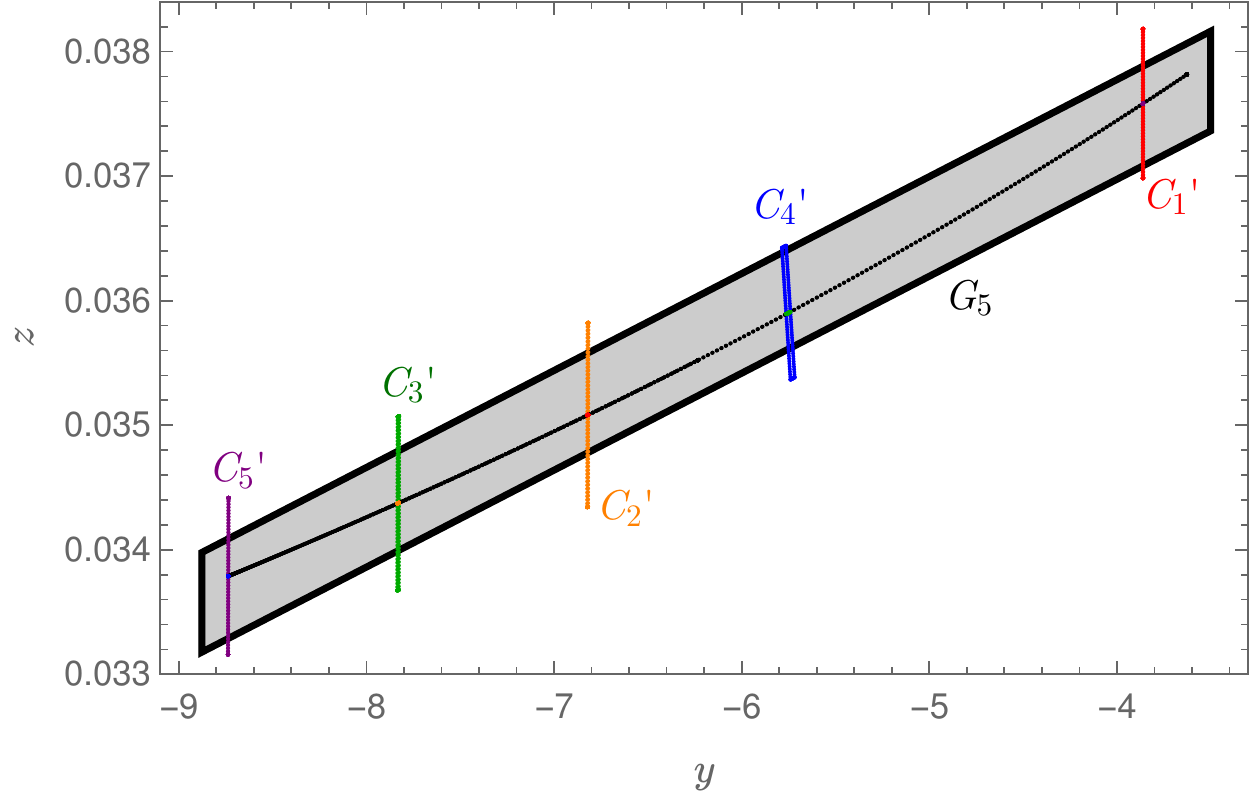}
	\caption{\label{fig:grid5}The contracting grid $G_5$ from Lemma \ref{lem:Roessler47}
	and its image through $P$. The supersets of $C_i$, $i=1,\dots,5$ and their images are marked in red, orange, green, blue and purple.}
\end{figure}

Let $a=4.7$. Then the system \eqref{eq:rossler} has an attracting $5$-periodic point for $P$ \cite[Lemma 7]{AGPZrossler}.
\begin{lemma}\label{lem:Roessler47}
The parallelogram $G_5$ in the $(y, z)$ coordinates on the section $\Pi$ (see
Fig.\ \ref{fig:grid5}):
\[
G_5= \bmatrix -6.1885 \\ 0.0356707 \endbmatrix
+ \bmatrix -1. & 0.000778356 \\ -0.000778356 & -1. \endbmatrix \cdot
\bmatrix \pm 2.68797 \\ \pm 0.0004 \endbmatrix
\]
is a contracting grid for the attracting $5$-periodic orbit from \cite[Lemma 7]{AGPZrossler}, with the inner cubes $C_i$, $i=1,\dots,5$, defined by $C_i = C'_i \cap G_5$, where
\begin{align*}
	&C'_1= \bmatrix -3.86108 \\ 0.0375827 \endbmatrix
+ \bmatrix 0.0693366 & 1. \\ -0.997593 & 0.000984231 \endbmatrix \cdot
\bmatrix \pm 0.0006 \\ \pm 0.00138 \endbmatrix
	\text{,}
	\\
	&C'_2= \bmatrix -6.82009 \\ 0.0350822 \endbmatrix
+ \bmatrix 0.7879108 & 1. \\ 0.615789 & 0.0007307 \endbmatrix \cdot
\bmatrix \pm 0.0012 \\ \pm 0.0024 \endbmatrix
	\text{,}
	\\
	&C'_3= \bmatrix -7.83056 \\ 0.0343732 \endbmatrix
+ \bmatrix 0.8138516 & 1. \\ 0.581073 & 0.000671 \endbmatrix \cdot
\bmatrix \pm 0.0012 \\ \pm 0.0042 \endbmatrix
	\text{,}
	\\
	&C'_4= \bmatrix -5.75153 \\ 0.0359038 \endbmatrix
+ \bmatrix 0.9997319 & 1. \\ -0.023153 & 0.0008062 \endbmatrix \cdot
\bmatrix \pm 0.0228 \\ \pm 0.01116 \endbmatrix
	\text{,}
	\\
	&C'_5= \bmatrix -8.73615 \\ 0.0337875 \endbmatrix
+ \bmatrix 0.8997843 & 1. \\ 0.436335 & 0.00062508 \endbmatrix \cdot
\bmatrix \pm 0.00144 \\ \pm 0.000744 \endbmatrix.
\end{align*}	
\end{lemma}

\begin{proof}
Computer-assisted, \cite{proof}.
\end{proof}

From the above Lemma and Theorem~\ref{th:Final} we obtain
\begin{theorem}(Compare to \cite[Theorem 6]{AGPZrossler})\label{th:r5}
The R\"ossler system \eqref{eq:rossler} with $a=4.7$ has $m$-periodic orbits for any $m\in \mathbb{N}\setminus\{3\}$ in the set $G_5$ defined in Lemma \ref{lem:Roessler47}.
\end{theorem}

\begin{remark}
In fact, Theorem~\ref{th:r5} is only a part of \cite[Theorem 6]{AGPZrossler}, because it does not exclude the possibility that a $3$-periodic orbit exists in $G_5$.
\end{remark}

\subsection{Case $a=4.381$}\label{ss:rossler6}
Consider $a=4.381$. From  Fig. \ref{fig:bif} it is clear that the R\"ossler system has an attracting $6$-periodic orbit.
\begin{lemma}\label{lem:6per}
The Poincar\'e map $P$ of the system \eqref{eq:rossler} with $a=4.381$ on the section $\Pi$ has a 6-periodic orbit $\mathcal{O}^6=\{p_1^6,p_2^6,p_3^6,p_4^6,p_5^6,p_6^6\}$, contained in the following rectangles in the $(y,z)$ coordinates on the section $\Pi$:
\begin{equation}
\begin{aligned}\label{eq:6-per}
		p_1^6 \in & -7.448265140_{33532}^{244187} \times 0.03638524011_{881493}^{973746}, \\	
		p_2^6 \in & -5.432682771_{276081}^{080253} \times 0.0381210024_{7833106}^{8150609}, \\
		p_3^6 \in & -8.146150765_{219835}^{118602} \times 0.03585533157_{361669}^{606319}, \\	
		p_4^6 \in & -4.05248247_{1003891}^{0816507} \times 0.03953831884_{313481}^{723778},\\	
		p_5^6 \in & -6.98865159_{7169091}^{6889441} \times 0.03675237717_{289087}^{31595},\\
		p_6^6 \in & -6.38538092_{5198882}^{4637889} \times 0.0372584624_{5305077}^{617641}.
\end{aligned}
\end{equation}
\end{lemma}

\begin{proof}
Computer-assisted \cite{proof}, via Interval Newton Method \cite{N}.
\end{proof}

\begin{figure}[h]
	\includegraphics[width=0.8\textwidth]{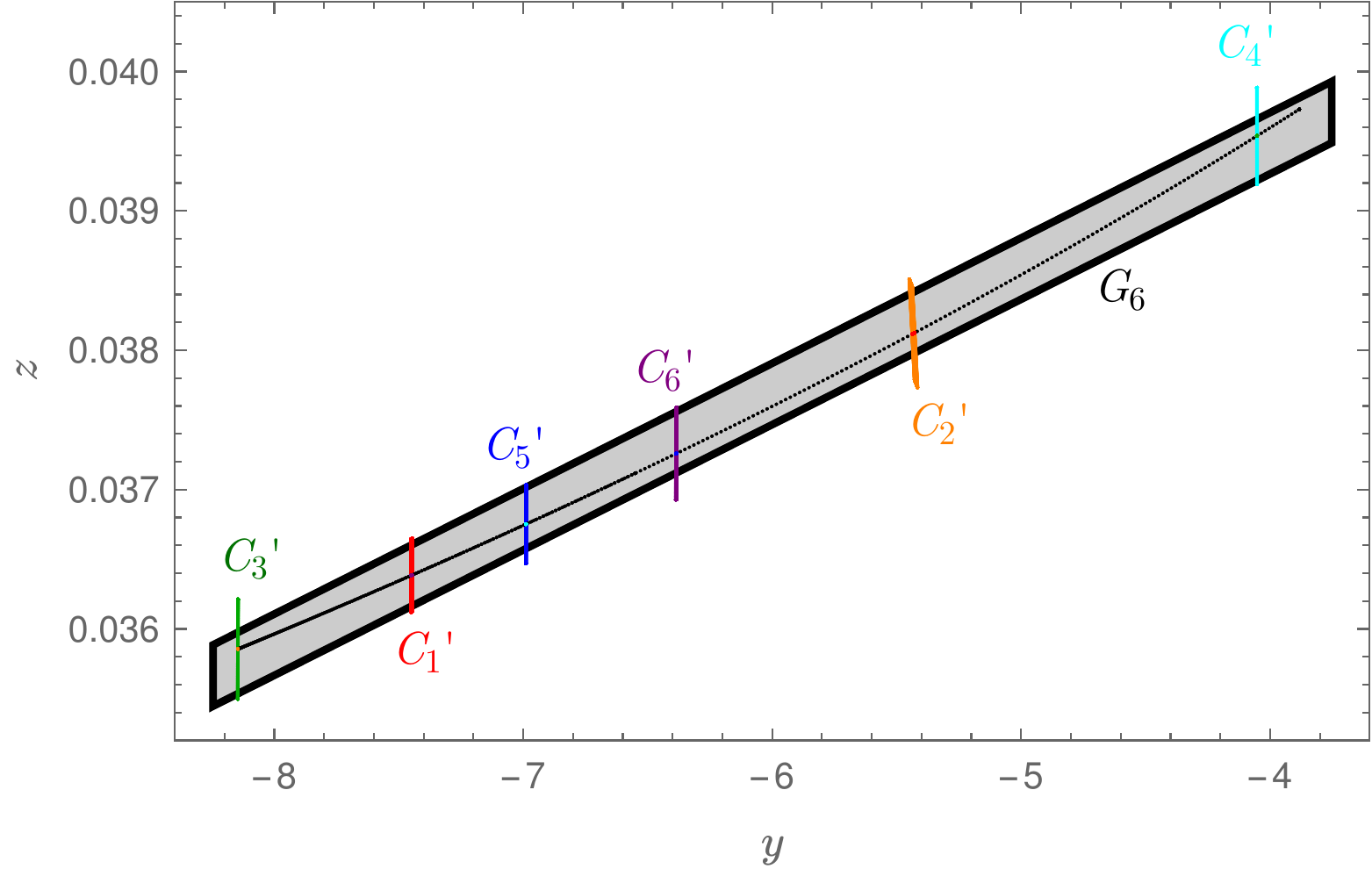}
	\caption{\label{fig:grid6}The contracting grid from Lemma \ref{lem:Roessler4381}
	and its image through $P$. The supersets of $C_i$, $i=1,\dots,6$ and their images are marked in red, orange, green, cyan, blue and purple.}
\end{figure}

\begin{lemma}\label{lem:Roessler4381}
The parallelogram $G_6$ in the $(y, z)$ coordinates on the section $\Pi$ (see
Fig.\ \ref{fig:grid6}):
\[
G_6= \bmatrix -5.99932 \\ 0.0376868 \endbmatrix
+ \bmatrix -1. & -0.000899679 \\ 0.000899679 & -1. \endbmatrix \cdot
\bmatrix \pm 2.24683 \\ \pm 0.00022 \endbmatrix
\]
is a contracting grid for the attracting $6$-periodic orbit \eqref{eq:6-per} from Lemma \ref{lem:6per}, with the inner cubes $C_i$, $i=1,\dots,6$, defined by $C_i = C'_i \cap G_6$, where
\begin{align*}
	&C'_1= \bmatrix -7.44827 \\ 0.0363852 \endbmatrix
+ \bmatrix1. & 0.8498 \\ 0.0007825 & 0.527106 \endbmatrix \cdot
\bmatrix \pm 0.00225 \\ \pm 0.0005 \endbmatrix
	\text{,}
	\\
	&C'_2= \bmatrix -5.43268 \\ 0.038121 \endbmatrix
+ \bmatrix 1.23042 & 0.696746 \\ -0.00567154 & -0.0240555 \endbmatrix \cdot
\bmatrix \pm 0.00509 \\ \pm 0.015 \endbmatrix
	\text{,}
	\\
	&C'_3= \bmatrix -8.14614 \\ 0.0358553 \endbmatrix
+ \bmatrix 1. & 0.907289 \\ 0.000736978 & 0.420507 \endbmatrix \cdot
\bmatrix \pm 0.000265 \\ \pm 0.00085 \endbmatrix
	\text{,}
	\\
	&C'_4= \bmatrix -4.05249 \\ 0.0395383 \endbmatrix
+ \bmatrix 0.999999 & 0.155044 \\ 0.00111181 & -0.987908 \endbmatrix \cdot
\bmatrix \pm 0.000485 \\ \pm 0.00035 \endbmatrix
	\text{,}
	\\
	&C'_5= \bmatrix -6.98865 \\ 0.0367524 \endbmatrix
+ \bmatrix 1. & 0.827066 \\ 0.000815525 & 0.562105 \endbmatrix \cdot
\bmatrix \pm 0.000712 \\ \pm 0.0005 \endbmatrix
	\text{,}
	\\
	&C'_6= \bmatrix -6.38538 \\ 0.0372585 \endbmatrix
+ \bmatrix 1. & 0.834783 \\ 0.000863145 & 0.550579 \endbmatrix \cdot
\bmatrix \pm 0.00149 \\ \pm 0.0006 \endbmatrix.
\end{align*}	
\end{lemma}

\begin{proof}
Computer-assisted, \cite{proof}.
\end{proof}

From the above Lemma and Theorem~\ref{th:Final} we obtain
\begin{theorem}\label{th:r6}
The R\"ossler system \eqref{eq:rossler} with $a=4.381$ has $m$-periodic orbits for any even $m\in 2\mathbb{N}$ and $m=1$, in the set $G_6$ defined in Lemma \ref{lem:Roessler4381}.
\end{theorem}

\begin{remark}
The above result cannot be strengthened to include odd periods because the permutation of the orbit with respect to the location on the model 1-dimensional manifold is as follows:
\[
\begin{tikzpicture}
\foreach \x/\l in {1/4, 2/2, 3/6, 4/5, 5/1, 6/3}
{
\coordinate (p\l) at (\x,0);
\coordinate (pg\l) at (\x,0.2);
\coordinate (pd\l) at (\x,-0.13);
\draw (\x,0.1) -- (\x,-0.1) node[below]{$p_\l$};
}
\draw[->] (0,0) -- (7,0);
\begin{scope}[>={stealth}]
\foreach \b/\e in {4/5, 2/3, 6/1}
\draw (pg\b) edge[->,bend left=40] (pg\e);
\foreach \b/\e in {5/6, 1/2, 3/4}
\draw (pd\b) edge[->,bend left=40] (pd\e);
\end{scope}
\end{tikzpicture}
\]
\end{remark}
Note that this is the case of an even-periodic orbit for which a \u{S}tefan sequence cannot be constructed (`all points switch sides', see details in \cite{Burns}) and using the method from the proof of Burns and Hasselblatt we obtain non-repeating loops of covering relations only of even length or a self-covering.

However, the next case shows the situation, in which one can prove the existence of odd-periodic points for an even-periodic attracting orbit. As we have mentioned, the difference lies in the permutation of the orbit.

\subsection{Case $a=5.42$}\label{ss:rossler6_2}
Let now see another case with an attracting $6$-periodic orbit for the R\"ossler system. Consider $a=5.42$ (see Fig. \ref{fig:bif}).
\begin{lemma}\label{lem:6per2}
The Poincar\'e map $P$ of the system \eqref{eq:rossler} with $a=4.381$ on the section $\Pi$ has a 6-periodic orbit $\mathcal{O}^6=\{p_1^6,p_2^6,p_3^6,p_4^6,p_5^6,p_6^6\}$, contained in the following rectangles in the $(y,z)$ coordinates on the section $\Pi$:
\begin{equation}
\begin{aligned}\label{eq:6-per2}
		p_1^6 \in & -3.3303887279_{60296}^{4934} \times 0.033810081022_{70888}^{86536}, \\	
		p_2^6 \in & -6.0438781482_{33535}^{13811} \times 0.0319883054102_{6752}^{8062}, \\
		p_3^6 \in & -9.930004688_{71182}^{693574} \times 0.029985122265_{72283}^{83182}, \\	
		p_4^6 \in & -3.561109751_{505439}^{469876} \times 0.033636111425_{11445}^{66204},\\	
		p_5^6 \in & -6.450138010_{324274}^{261234} \times 0.0317509776490_{3493}^{7362},\\
		p_6^6 \in & -10.061811798_{91221}^{88145} \times 0.02992604451_{798922}^{853264}.
\end{aligned}
\end{equation}
\end{lemma}

\begin{proof}
Computer-assisted \cite{proof}, via Interval Newton Method \cite{N}.
\end{proof}

\begin{figure}[h]
	\includegraphics[width=0.8\textwidth]{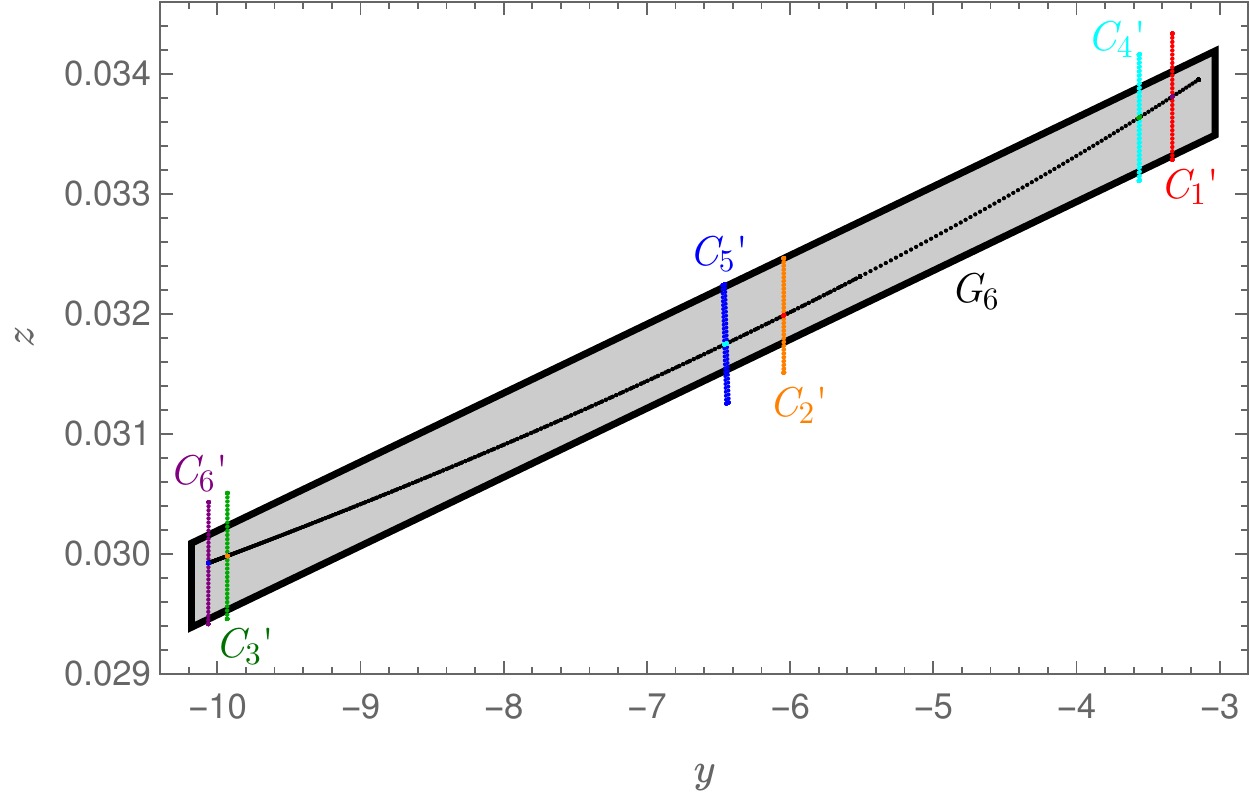}
	\caption{\label{fig:grid6_2}The contracting grid from Lemma \ref{lem:Roessler542}
	and its image through $P$. The supersets of $C_i$, $i=1,\dots,6$ and their images are marked in red, orange, green, cyan, blue and purple.}
\end{figure}

\begin{lemma}\label{lem:Roessler542}
The parallelogram $G_6$ in the $(y, z)$ coordinates on the section $\Pi$ (see
Fig.\ \ref{fig:grid6_2}):
\[
G_6= \bmatrix -6.60556 \\ 0.0317909 \endbmatrix
+ \bmatrix -1. & -0.000573253 \\0.000573253 & -1. \endbmatrix \cdot
\bmatrix \pm 3.57445 \\ \pm 0.00035 \endbmatrix
\]
is a contracting grid for the attracting $6$-periodic orbit \eqref{eq:6-per2} from Lemma \ref{lem:6per2}, with the inner cubes $C_i$, $i=1,\dots,6$, defined by $C_i = C'_i \cap G_6$, where
\begin{align*}
	&C'_1= \bmatrix -3.33039 \\ 0.0338101 \endbmatrix
+ \bmatrix 1. & 0.0114844 \\ 0.000763188 & -0.999934 \endbmatrix \cdot
\bmatrix \pm 0.0015225 \\ \pm 0.000525 \endbmatrix
	\text{,}
	\\
	&C'_2= \bmatrix -6.04388 \\ 0.0319883 \endbmatrix
+ \bmatrix 1. & 0.566012 \\ 0.000593828 & -0.824397 \endbmatrix \cdot
\bmatrix \pm 0.0029925 \\ \pm 0.0005775 \endbmatrix
	\text{,}
	\\
	&C'_3= \bmatrix -9.93 \\ 0.0299851 \endbmatrix
+ \bmatrix 1. & 0.866643 \\ 0.000450065 & 0.498928 \endbmatrix \cdot
\bmatrix \pm 0.0021 \\ \pm 0.00105 \endbmatrix
	\text{,}
	\\
	&C'_4= \bmatrix -3.56111 \\ 0.0336361 \endbmatrix
+ \bmatrix 1. & 0.0148011 \\ 0.000745026 & -0.99989 \endbmatrix \cdot
\bmatrix \pm 0.0043575 \\ \pm 0.000525 \endbmatrix
	\text{,}
	\\
	&C'_5= \bmatrix -6.45014 \\ 0.031751 \endbmatrix
+ \bmatrix 1. & 0.999296 \\ 0.000574732 & -0.0375247 \endbmatrix \cdot
\bmatrix \pm 0.00945 \\ \pm 0.013125 \endbmatrix
	\text{,}
	\\
	&C'_6= \bmatrix -10.0618 \\ 0.029926 \endbmatrix
+ \bmatrix 1. & 0.887687 \\ 0.000446372 & 0.460448 \endbmatrix \cdot
\bmatrix \pm 0.00084 \\ \pm 0.0011025 \endbmatrix.
\end{align*}	
\end{lemma}

\begin{proof}
Computer-assisted, \cite{proof}.
\end{proof}

From the above Lemma and Theorem~\ref{th:Final} we obtain
\begin{theorem}
The R\"ossler system \eqref{eq:rossler} with $a=5.42$ has $m$-periodic orbits for any even $m\in 2\mathbb{N}$ and $m=1$, in the set $G_6$ defined in Lemma \ref{lem:Roessler542}.
\end{theorem}

\begin{remark}
In this case, we are able to strengthen the above result. Note that the permutation of the orbit with respect to the location on the model 1-dimensional manifold is the following:
\begin{equation}\label{eq:6-permutation}
\begin{tikzpicture}
\foreach \x/\l in {1/6, 2/3, 4/5, 5/2, 7/4, 8/1}
{
\coordinate (p\l) at (\x,0);
\coordinate (pg\l) at (\x,0.2);
\coordinate (pd\l) at (\x,-0.13);
\draw (\x,0.1) -- (\x,-0.1) node[below]{$p_\l$};
}
\foreach \b/\e/\l in {4/2/0, 5/3/1}
\draw[ultra thick] (p\b) -- (p\e) node[midway,above]{$I_\l$};
\draw[->] (0,0) -- (9,0);
\begin{scope}[>={stealth}]
\foreach \b/\e in {1/2, 2/3, 4/5, 5/6}
\draw (pg\b) edge[->,bend right=40] (pg\e);
\foreach \b/\e in {3/4, 6/1}
\draw (pd\b) edge[->,bend right] (pd\e);
\end{scope}
\end{tikzpicture}
\end{equation}
and, although the orbit is even-periodic, we are able to construct a \u{S}tefan sequence. One of the shortest possible diagrams of 1-dimensional covering relations that we obtain involves only two $\mathcal{O}$-intervals: $I_0 = [p_2,p_4]$ and $I_1 = [p_3,p_5]$, as denoted on \eqref{eq:6-permutation}. Then the diagram is
	\begin{equation}\label{eq:diag2}
	\begin{tikzcd}
I_0
	\arrow[r, rightarrow, bend left=15]
&
I_1
	\arrow[l, rightarrow, bend left=15]
	\arrow[rightarrow, loop, distance=2em, in=-20, out=20]
\end{tikzcd}
\text{,}
\end{equation}	
and all covering relations are proper. Therefore, using the grid from Lemma \ref{lem:Roessler542}, we can construct a similar diagram of horizontal 2-dimensional covering relations:
\begin{equation}
\begin{tikzcd}
S(I_0)
	\arrow[r, Rightarrow, bend left=10]
&
S(I_1)
	\arrow[l, Rightarrow, bend left=10]
	\arrow[Rightarrow, loop, distance=2em, in=-10, out=15]
\end{tikzcd}
\text{,}
\end{equation}
from which follows the existence of $m$-periodic points for $P$, for all natural $m$.
\end{remark}

Therefore we obtain
\begin{theorem}
The R\"ossler system \eqref{eq:rossler} with $a=5.42$ has $m$-periodic orbits for any $m\in \mathbb{N}$ in the set $G_6$ defined in Lemma \ref{lem:Roessler542}.
\end{theorem}

\section*{Acknowledgment}
The second author is supported by the Funding Grant NCN UMO-2016/22/A/ST1/00077.


\bibliographystyle{plain}
\bibliography{periodic_bib}

\begin{thebibliography}{10}

\bibitem{Block}
L.~Block, J.~Guckenheimer, M.~Misiurewicz, and L.~Young.
\newblock Periodic points and topological entropy of one dimensional maps.
\newblock {\em Lecture Notes Math.}, 819:18, 11 2006.

\bibitem{Burns}
Keith Burns and Boris Hasselblatt.
\newblock The {S}harkovsky theorem: A natural direct proof.
\newblock {\em The American Mathematical Monthly}, 118(3):229--244, 2011.

\bibitem{capd}
{CAPD group}.
\newblock {\em Computer Assisted Proofs in Dynamics C++ library}.
\newblock \url{http://capd.ii.uj.edu.pl}.

\bibitem{proof}
A.~Gierzkiewicz and P.~Zgliczy\'nski.
\newblock {C}++ source code.
\newblock Available online on \\
  \url{http://kzm.ur.krakow.pl/~agierzkiewicz/publikacje.html}.

\bibitem{AGPZrossler}
A.~Gierzkiewicz and P.~Zgliczy{\'{n}}ski.
\newblock Periodic orbits in the {R}\"ossler system.
\newblock {\em Communications in Nonlinear Science and Numerical Simulation},
  page 105891, 2021.

\bibitem{capd-article}
T.~Kapela, M.~Mrozek, D.~Wilczak, and P.~Zgliczy\'nski.
\newblock \texttt{CAPD::DynSys}: a flexible \texttt{C++} toolbox for rigorous
  numerical analysis of dynamical systems.
\newblock {\em Communications in Nonlinear Science and Numerical Simulation},
  page 105578, 2020.

\bibitem{N}
A.~Neumaier.
\newblock {\em Interval Methods for Systems of Equations}.
\newblock Encyclopedia of Mathematics and its Applications. Cambridge
  University Press, 1991.

\bibitem{Rossler76}
O.E. R\"ossler.
\newblock An equation for continuous chaos.
\newblock {\em Physics Letters A}, 57(5):397 -- 398, 1976.

\bibitem{ShU}
A.N. Sharkovskii.
\newblock Co-existence of cycles of a continuous mapping of the line into
  itself.
\newblock {\em Ukrainian Math. J.}, 16:61--71, 1964.
\newblock (in Russian, English translation in J. Bifur. Chaos Appl. Sci.
  Engrg., 5:1263-1273, 1995.).

\bibitem{Stefan}
P.~\u{S}tefan.
\newblock A theorem of {S}arkovskii on the existence of periodic orbits of
  continuous endomorphisms of the real line.
\newblock {\em Comm. Math. Phys.}, 54(3):237--248, 1977.

\bibitem{PZmulti}
P.~Zgliczy\'nski.
\newblock Multidimensional perturbations of one-dimensional maps and stability
  of \u{S}arkovsk\u{i} ordering.
\newblock {\em International Journal of Bifurcation and Chaos},
  09(09):1867--1876, 1999.

\bibitem{PZszarI}
P.~Zgliczy\'nski.
\newblock Sharkovskii's theorem for multidimensional perturbations of
  one-dimensional maps.
\newblock {\em Ergodic Theory and Dynamical Systems}, 19(6):1655--1684, 1999.

\bibitem{GZ}
P.~Zgliczy\'nski and M.~Gidea.
\newblock Covering relations for multidimensional dynamical systems.
\newblock {\em Journal of Differential Equations}, 202(1):32--58, 2004.

\end{thebibliography}

\end{document}